\input amstex\documentstyle{amsppt}  
\pagewidth{12.5cm}\pageheight{19cm}\magnification\magstep1
\topmatter
\title History of the canonical
basis and crystal basis
\endtitle
\author George Lusztig\endauthor
\endtopmatter   
\document

\define\part{\partial}

\define\iy{\infty}

\define\ti{\tilde}

\redefine\i{^{-1}}

\redefine\l{\lambda}

\define\BB{\bold B}

\define\LL{\bold L}

\define\NN{\bold N}

\define\QQ{\bold Q}

\define\UU{\bold U}

\define\ZZ{\bold Z}

\define\ca{\Cal A}

\define\cl{\Cal L}

\define\tb{\ti b}

\head 1. Introduction\endhead
This document clarifies the historical
development of the canonical basis theory (see \S3).
Such a clarification is necessary in view of the
persistent false and misleading claims in the literature with regard to the history of the
canonical basis and crystal basis
(see \S2).
Over the years I thought that the evidence in the literature
about this history is clear and that most people accept this
evidence. However Kashiwara repeatedly made misleading
claims about priority and this has influenced some people who do not know the history.
Since such false and  misleading information has now
spread out of the mathematical community to the public, I feel that I must now publish
this document. I present this document to restore the
correct historical record and to defend academic
integrity in the mathematical community.

\head 2. False claims\endhead
\subhead 2.1. An evolving narrative \endsubhead

In 1995, 
 Kashiwara \cite{K95} claimed  that his basis and
mine were ``constructed independently'' despite clear
evidence to the contrary, see \S3.

In 2006 at my birthday conference he gave a talk where
he again made the independence claim.
After his talk I pointed out this error
to him and then
he \bf acknowleged in \cite{EK} that I defined
the canonical basis first for type ADE.\rm

Yet, in 2018, Kashiwara, \cite{K18} (version 1 on arXiv),
repeated the independence claim. After 
 I protested, he \bf acknowleged in version 2 on
 arXiv that I defined the canonical basis first for
 type ADE.\rm

Kashiwara's repeated attempts to rewrite  history
have also influenced other people who do not know the history.
See 2.2 and 2.3 for examples of this.

\subhead 2.2. The Carter incident\endsubhead
In the original
version of \cite{C}, Carter states: ``In 1990 Lusztig made another
discovery of fundamental importance by proving the existence of a
remarkable basis of a quantized enveloping algebra called the
canonical basis. This basis was also subsequently proved to exist
by Kashiwara...''
However, in the published version,
\bf somehow
the word ``subsequently'' was changed to
``independently'', without Carter's
consent, \rm falsely implying parallel discovery.
After I pointed out this distorsion to the Editorial
Committee, \bf they issued an
Erratum \cite{EC},
restoring Carter's original wording.\rm

\subhead 2.3. Misleading awards citations\endsubhead
Throughout the past eight years, Kashiwara
has received various awards based in part
on inaccurate information provided by nominators.

-In the Chern Medal initial citation it is stated that
``Kashiwara's discovery of crystal basis is
another landmark in representation theory''.

-In the Kyoto Prize initial citation it is stated that
``...Kashiwara...for construction of the crystal basis theory''.

-In the Abel Prize initial citation it is stated that
``Kashiwara introduced the notion of crystal bases and
proved the existence of crystal bases for integrable highest weight modules for quantum groups...Kashiwara also generalized crystal bases to global bases which were independently discovered by George Lusztig under the name
canonical bases''.

\bf Each of these statements is false, see 3.5, 3.6. 
In each of these cases the citations were modified after I
explained the facts. \rm

\subhead 2.4. Public misinformed\endsubhead
In the New York Times (March 26, 2025) it is stated that Kashiwara
invented the crystal basis.

In Nature (March 26, 2025) it is stated: ``In particular,
Kashiwara's
notion of a crystal base has enabled mathematicians to interpret any representation ...''

In Scientific American
(March 26, 2025) it is stated:``...Kashiwara introduced
the concept of crystal bases.''

In the Wikipedia page on Abel Prize it is stated:
``(Kashiwara) ... and the discovery of crystal bases''.

These statements reflect the statements
in the original Abel Prize citation and that
the public was misinformed.

\head 3. Overview of the canonical basis\endhead
\subhead 3.1\endsubhead
In my 1990 paper \cite{L90} I observed that for the
irreducible representations of a quantized enveloping algebra of
simply laced type as well as for the positive part of that
algebra there is a new, extremely rigid structure in
which the objects of the theory are provided with canonical bases
with rather remarkable properties; in particular, all the structure
constants with respect to the canonical basis are in
$\bold N[v,v^{-1}]$ ($v$ is an indeterminate).

My construction employed two methods: an algebraic approach utilizing
braid group actions and PBW (Poincar\'e-Birkhoff-Witt) bases and a
topological approach using intersection cohomology.

Key consequences include:

Specializing the parameter $v$ to $v=1$ yields canonical bases for the irreducible representations of the corresponding simple Lie algebra
in which all the structure constants are in $\bold N$,
a property that was not previously known. When specializing
$v$ to $0$ or to $\iy$ the canonical basis produces
a shadow structure (``crystal basis'') which Kashiwara also considered
(for classical types) in his 1990 paper \cite{K90}.

The graph structure on the crystal basis considered by Kashiwara
is also a shadow of the canonical basis, see \cite{L90b, Theorem 7.5}, in the sense that the graph structure
can be deduced from the canonical basis. A property
similar to \cite{L90b, Theorem 7.5} is stated in
Kashiwara's later paper \cite{K95}, Lemma 12.1, but he fails to give a reference to my paper.

\subhead 3.2\endsubhead
In \cite{K90}, written at the same time
as \cite{L90}, Kashiwara gives a conjectural definition of the crystal basis
(``basis at $v=0$'') for the irreducible
representations of a quantized enveloping algebra. He proves that the definition is correct for classical types. He was motivated by mathematical physics and his focus on looking at $v=0$ is related to the absolute temperature being $0$. The underlying philosophy,
as indicated in the paper \cite{K90}, is that at the absolute
temperature being $0$, modules should have some special properties. This philosophy would not lead to canonical bases.

\subhead 3.3\endsubhead
His subsequent announcement \cite{K90a} -which
explicitly cites \cite{L90} and states that his basis
``appears to agree'' with mine-proves that he
encountered the concept in my work first.
He developed his approach after reading \cite{L90}
though without crediting this foundation.

The main ingredients in lifting the crystal basis to the canonical
bases in \cite{L90} are

-the $\ZZ[v,v\i]$-form of quantized enveloping algebras
\cite{L88, L90a},

-the bar involution \cite{L90}.

\bf In \cite{K90a} Kashiwara copies (see  (1.5), page 278 and \S5,
page 279) \rm
\bf the definitions of these ingredients from
\cite{L88}, \cite{L90}, \cite{L90a}  without reference, \rm
\bf thus giving the impression that he is the first to define these concepts.\rm

\bf He also copies (see the line after Theorem 5, page 279) the method of \rm
\bf lifting the crystal basis to the canonical basis from \cite{L90} without reference.\rm

\subhead 3.4\endsubhead
The concept of canonical basis was subsequently extended to the
broader setting of Kac-Moody Lie algebras. I provided such an extension
using a topological approach in \cite{L91} (see also \cite{L93})
while Kashiwara extended it
through an algebraic approach in \cite{K91}.

\subhead 3.5\endsubhead
It is \bf not \rm
correct to say that Kashiwara defined
the canonical basis independently of me. 
The reason is as follows.

My 1990 paper \cite{L90} contains the first construction of
canonical bases for $ADE$ types.

Kashiwara's 1990 paper \cite{K90} written concurrently
with \cite{L90} contains nothing about canonical
bases. 
Both the motivation and the technical tools of his paper \cite{K91} are inherited
from \cite{L90}, although Kashiwara does
not acknowledge that.

Therefore it is clear that without my discovery of canonical bases \cite{L90},
Kashiwara could not have written \cite{K91}.

\subhead 3.6\endsubhead
It is \bf not \rm
correct to say that Kashiwara is the sole
discoverer of the crystal basis. The reason is
as follows.

My 1990 paper \cite{L90} defines the crystal basis
for $ADE$ types, including the most challenging $E_8$ type.

Kashiwara's 1990 paper \cite{K90}
defines the crystal basis for classical types but
in type E his definition is only conjectural.

In \cite{K91}, \cite{L91} (see also \cite{L93})
the crystal basis is defined for all types.
But to do so, \cite{K91} needed the idea that the crystal basis
exists for the $+$ part of the quantum group (not only for
irreducible representations, which was the only case considered
in \cite{K90}). This idea came from my paper
\cite{L90},
without acknowledgment.

\head 5. Conclusion\endhead
By publishing this document I aim to rectify
the historical narrative for the benefit of the mathematical community and of the general
public and to ensure that proper attribution and academic integrity is upheld by all.

I trust that all readers -including Kashiwara- 
will recognize these established facts:

(a) The canonical basis was first defined in my work
\cite{L90} and Kashiwara's subsequent contribution
built directly on this foundation.

(b) The crystal basis is not solely Kashiwara's discovery.

And everyone who knows the history would suggest Kashiwara to publicly acknowledge (a) and (b), to correct all false and misleading
information once and for all.

\head Appendix.
Technical details and comparisons\endhead
\subhead A.1\endsubhead
In this appendix we give more details of the papers
\cite{L90}, \cite{K90}, \cite{K90a}.

We introduce some notation.
Let $\UU$ be the quantum group over $\QQ(v)$
attached by Drinfeld, Jimbo to a generalized
Cartan matrix $C=(a_{ij})_{i,j\in I}$.
Let $E_i,F_i,K_i$ ($i\in I$) be the standard generators of $\UU$;
recall that $\UU^+$ is the subalgebra of $\UU$ generated by
$E_i,i\in I$.
Now let $\l:I@>>>\NN$ and let
$\LL_\l=\UU^+/\sum_{i\in I}\UU^+E_i^{\l(i)+1}$. Let $x_0$ be the
image of $1\in\UU^+$ in $\LL_\l$. It is known
that $\LL_\l$ is an irreducible $\UU$-module in which $E_i$ acts by
left multiplication, $F_i$ maps $x_0$ to $0$ and $K_i$ maps $x_0$ to
$v^{-\l(i)}x_0$.

There is a well defined $\QQ$-algebra isomorphism (``bar involution'')
$\bar{}:\UU^+@>>>\UU^+$ such that $\bar E_i=E_i,\bar v=v\i$.

Let $\ca'=\QQ[v,v\i]$.
In \cite{L88} I defined an $\ca'$-form  $\UU^+_{\ca'}$ of $\UU^+$.

\subhead A.2\endsubhead
We now assume that $C$ is symmetric, positive definite.
Let $\ca=\ZZ[v,v\i]$. In \cite{L90a}, I defined an $\ca$-submodule
$\UU^+_\ca$ of $\UU^+$ which was an $\ca$-subalgebra and defined
an $\ca$-basis for it. This is one of several $\ca$-bases
(PBW bases) which can be defined using the braid group action
\cite{L88} on $\UU$ (there is one PBW basis for each reduced
expression of the longest element of the Weyl group).

\subhead A.3\endsubhead
The paper \cite{L90} contains the definition of the canonical basis
of $\UU^+$. The definition is as follows.

We show  that the $\ZZ[v\i]$-submodule $\cl$ of $\UU^+$
generated by any PBW basis is independent of the PBW basis and that
the image of any
PBW basis under the obvious map $\pi:\cl@>>>\cl/v\i\cl$ is a
$\ZZ$-basis $\beta$ of $\cl/v\i\cl$ independent of the PBW-basis.

We show that

(a) for any $b\in\beta$ there is a unique
$\tb\in\cl$ such that $\bar{\tb}=b,\pi(\tb)=b$;

(b) $\BB:=\{\tb;b\in\beta\}$ is a
$\QQ(v)$-basis of $\UU^-$, an $\ca$-basis of $\UU^+_\ca$,
a $\ZZ[v\i]$-basis of $\cl$ and a $\ZZ$-basis of $\cl\cap\bar\cl$.

This is the canonical basis of $\UU^+$.

\subhead A.4\endsubhead
Let $\l:I@>>>\NN$.  In \cite{L90} it is shown that
the nonzero vectors
in the image of $\BB$ under the obvious map $\UU^+@>>>\LL_\l$ form a
basis $\BB_\l$ of $\LL_\l$. This is the canonical basis of $\LL_\l$.
By specializing $v$ to $1$ one obtains a canonical basis of any
irreducible representation of the Lie algebra defined by $C$.

\subhead A.5\endsubhead
We return to a general $C$. Let $\l:I@>>>N$.
Let $A$ be the ring consisting of the
elements in $\QQ(v)$ which have no pole at
$v=0$. In \cite{K90} an explicit collection $X_\l$ of vectors in
$\LL_\l$ is defined. Let $L(\l)$ be the $A$-submodule of $\LL_\l$
spanned by $X_\l$ and let $B(\l)$
be the set of nonzero elements in the
image of $X_\l$ under the obvious map $L(\l)@>>>L(\l)/vL(\l)$.
In \cite{K90} it is conjectured that $B(\l)$ is a $\QQ$-basis of
$L(\l)/vL(\l)$. This conjecture is proved in \cite{K90} in the case
where $C$ is of finite, classical type.

\subhead A.6\endsubhead
In \cite{K90a}, Kashiwara 
announced a proof of his conjecture in A.5 and a proof of the following
version of that conjecture.

(a) There is an explicit collection $X_\iy$ of vectors in
$\UU^+$ such that if $L(\iy)$ is the $A$-submodule of $\UU^+$
spanned by $X_\iy$ and $B(\iy)$
is the set of nonzero elements in the
image of $X_\iy$ under the obvious map
$\pi_\iy:L(\iy)@>>>L(\iy)/vL(\iy)$
then $B(\iy)$ is a $\QQ$-basis of $L(\iy)/vL(\iy)$.

Now $L(\iy),B(\iy)$ did not appear in
\cite{K90}. They did appear in \cite{L90} for type ADE; indeed
in that case $L(\iy)$ is the same as $A\otimes\cl$
(with $\cl$ as in A.3, after changing $v$ to $v\i$) and $B(\iy)$
is similarly the same as $\beta$ in A.3. But Kashiwara
does not say that the idea to consider $L(\iy),B(\iy)$ (which was
missing in \cite{K90}) comes from \cite{L90}.

\subhead A.7\endsubhead
The paper \cite{K90a} announces a proof of 
statements (a),(b) below.

(a) For any $b\in B(\iy)$ there is a unique
$\tb\in L(\iy)\cap\UU^-_{\ca'}$ such that
$\bar{\tb}=\tb,\pi_\iy(\tb)=b$.

(b) The set $\BB'=\{\tb;b\in B(\iy)\}$ is
a $\QQ(v)$-basis of $\UU^+$ and an
$\ca'$-basis of $\UU^+_{\ca'}$.

$\BB'$ is called the global crystal basis of $\UU^+$.
There is also a statement about a
global crystal basis $\BB'_\l$ of $\LL_\l$ for $\l:I@>>>\NN$.

\subhead A.8\endsubhead
The
idea to lift $b$ to $\tb$ as in A.7(a) is copied in
\cite{K90a}
from the analogous idea in \cite{L90}, see A.3(a),
without reference.

\subhead A.9\endsubhead
The statements of \cite{K90a} were proved in \cite{K91}.
A geometric construction of the canonical basis $\BB$
for symmetric $C$
(which has very strong positivity properties not seen in the
approach of \cite{K91}) is given in \cite{L91} extending one of the
two approaches of \cite{L90}. The case of
symmetrizable $C$ is deduced from the symmetric case in \cite{L93}.
It is known that $\BB=\BB'$ (see \cite{L90b} for the case where
$C$ is symmetric, positive definite and \cite{GL} for the other cases).

\enddocument